\documentclass[11pt]{article}
\usepackage{amssymb}
   \csname @addtoreset\endcsname{equation}{section}
\newtheorem{theorem}{Theorem}%[section]
\newtheorem{lemma}{Lemma}%[section]
\def\e{\varepsilon}

\def\defi{\stackrel{{\scriptscriptstyle \Delta}}{=}}

\def\N{{\cal N}}

\def\o{\omega}

\def\F{{\cal F}}
\def\w{\widehat}

\def\R{{\bf R}}

\def\X{{\cal X}}

\def\M{{\cal M}}

\newcommand{\be}{\begin{equation}}
\newcommand{\ee}{\end{equation}}
\newcommand{\bd}{\begin{displaymath}}
\newcommand{\ed}{\end{displaymath}}
\newcommand{\ba}{\begin{array}{ll}}
\newcommand{\ea}{\end{array}}
\newcommand{\baa}{\begin{eqnarray}}
\newcommand{\eaa}{\end{eqnarray}}
\newcommand{\baaa}{\begin{eqnarray*}}
\newcommand{\eaaa}{\end{eqnarray*}}
\font\sm=cmr10
\def\iomega{i\o}
\title{
On identity theorem for real functions %XXX
 }
\author{ Nikolai Dokuchaev\\ {\sm  Department of
Mathematics, Trent University, Ontario, Canada}}
 
\begin{document}
 \vspace{-0.5cm}
 \maketitle
\begin{abstract} Identity theorem for analytic  complex functions says that a function
is uniquely defined by its values on a set that contains  a density
point.   The paper presents sufficient conditions for classes of
real analytic functions  that ensures similar property.
\\ AMS 2010 classification :
26A99,    %   None of the above, but in this section
42A38 %Fourier and Fourier-Stieltjes transforms and other transforms
%of Fourier type
%\\
%PACS 2008 numbers: 02.30.Mv, %    Approximations and expansions
%02.30.Nw,  %  Fourier analysis
%02.30.Yy, %    Control theory
%07.05.Mh,  %  Neural networks, fuzzy logic, artificial intelligence
%07.05.Kf %Data analysis: algorithms and implementation; data management
\end{abstract}
\section{Introduction} Identity theorem for complex functions states that an analytic in a domain
function is uniquely defined by its values on a subdomain. It is
well known that the real functions do not have this feature. For
instance, let $x_1(t)= 0$ for all $t\in\R$ and let $x_2(t)$ be the
so-called Sobolev kernel defined as $x_2(t)=\exp(t^2(1-t^2)^{-1})$,
$t\in (-1,1)$ and $x_2(t)=0$ if $t\in\R\backslash(-1,1)$. Both
functions are infinitely differentiable, and $x_1(t)=x_2(t)$ for
$t\le -1$. This example shows that the identity property does not
hold for infinitely differentiable real functions. On the other
hand, there are some classes of real  functions that are uniquely
defined by their values on a part of the real axis (for instance,
periodic functions and band-limited functions).  It is interesting
to find wider classes of real functions where identity theorem
holds. It may be useful for the extrapolation problems and
forecasts.
\par
  The paper presents sufficient conditions for some classes of
real analytic functions  that ensures identity if the functions are
identical on a semi-infinite interval. These conditions are
expressed in the term of boundaries for growth for $L_2$-norms of
derivatives as well as some integrability condition for the Fourier
transforms.
\section{Problem setting and main result}
 For $C>0$, consider a class $\M(C)$ of infinitely differentiable functions  $x(t):\R\to\R$
 such that there exists $M=M(x(\cdot))>0$ such that \baaa
&&\left\|\frac{d^kx}{dt^k}(\cdot)\right\|^2_{L_2(\R)}\le C^kM, \quad
k=0,2,4,6,.... \label{deriv1} \eaaa \par For $C>0$, consider a class
$\N(C)$  of  functions $x(t):\R\to\R$
  such that there exists $M=M(x(\cdot))>0$ such that \baa
\frac{1}{2}\left(\left\|\frac{d^{k-1}x}{dt^{k-1}}(\cdot)\right\|^2_{L_2(\R)}
+\left\|\frac{d^{k+1}x}{dt^{k+1}}(\cdot)\right\|^2_{L_2(\R)}\right)\le
k!C^{-k}M, \quad k=1,3,5,7,....\hphantom{x} \label{deriv2}\eaa
\par
 For  $q\in\{1,2\}$, let $\X(q)=\X(q,T)$ be the  set of
processes $x(\cdot)\in L_2(\R)\cup L_1(\R)$ such that
\baaa\int_{-\infty}^{+\infty}e^{qT|\o|}|X(\iomega)|^q d\o <+\infty,
\eaaa where $X(\iomega)=\F x$ is the Fourier transform of
$x(\cdot)$.
\par
Let $\M\defi \cup_{C>0}\M(C)$,  $\N\defi\cup_{C>0}\N(C)$, and
$\X(q)\defi\cup_{T>0}\X(q,T)$.
\begin{theorem}\label{Th1} Let
$x(\cdot)\in\M\cup \N\cup \X(2)$ be such that $x(t)=0$ for $t<0$.
Then $x(t)\equiv 0$.
\end{theorem}
\par
{\it Proof of Theorem \ref{Th1}.} By Propositions 1 and 2 \cite{B},
$\M\subset \X(2)$ and $\N\subset \X(2)$. (In this proposition, the
set of infinitely differentiable functions is denoted as
$C^{\infty}(\R)$). Therefore, it suffices to prove the theorem  for
$x(\cdot)\in \X(2)$.
\par
Let $T>0$ be such that $x(\cdot)\in \X(2,T)$.
\par
The following lemma is a special case of Theorem 1 \cite{B}.
\begin{lemma}\label{lemma} The process $x(\cdot)\in\X(2,T)$ is  weakly predictable in the
following sense: for any $\e>0$ and any kernel $k\in L_\infty(0,T)$,
there exists a kernel $\w k(\cdot)\in L_2(0,+\infty)\cap
L_\infty(0,+\infty)$ such that
$$
\|y-\w y\|_{L_2(\R)}\le \e,
$$
where
 \baaa y(t)\defi
\int_t^{t+T}k(t-s)x(s)ds,\qquad \w y(t)\defi \int^t_{-\infty}\w
k(t-s)x(s)ds.\label{predict} \eaaa
\end{lemma}
\par
We will use this lemma to prove Theorem \ref{Th1}. First, let us
observe that \baa \w y(t)=0\quad \forall t<0. \label{zero} \eaa
\par
Further, let us show that $x(t)=0$ for $t\in [0,T]$.  Let
$\{k_i(\cdot)\}_{i=1}^{+\infty}$ be a basis in $L_2(-T,0)$, with
continuous bounded functions $k_i$. Let $y_i(t)\defi\int_t^{t+T}
k_i(t-s)x(s)ds$. By Lemma \ref{lemma}, it follows from (\ref{zero})
that $y_i(\cdot)|_{t\le 0}=0$ as an element of $L_2(-\infty,0)$.
Since $y_i(t)$ is a continuous function, it follows that $y_i(t)=0$
for $t\le 0$. It follows that $x(\cdot)|_{[0,T]}= 0$ as an element
of $L_2(0,T)$. By the properties of the class $\X(2)$,  it follows
that $x(t)$ is continuous. Hence $x(t)=0$ for $t\le T$.
\par
Further, let us apply the proof given above to the function
$x_1(t)=x(t+T)$. Clearly, $x_1(\cdot)\in\X(2)$. We found that
$x_1(t)=0$ for $t<0$. Similarly, we obtain that $x_1(t)=0$ for all
$t\le T$, i.e., $x(t)=0$ for all $t<2T$. Repeating this procedure
$n$ times, we obtain that i.e., $x(t)=0$ for all $t<nT$ for all
$n\ge 1$. This completes the proof of Theorem \ref{Th1}. $\Box$
\section{Concluding remarks}
We established identity theorem for the functions with exponential
rate of decay for the spectrum. We don't claim that this rate of
decay is a necessary condition for identity theorem. However, it can
be observed that the polynomial rate of decay for spectrum is
insufficient for identity property. For instance, consider a process
$x(\cdot)\in L_2(\R)$ with Fourier transform $X(\iomega)=
(1+i\o)^{-N-1}$, where $N>0$ is an arbitrarily large integer. Then
$x(t)=0$ for all $t\le 0$ and
$$
\int_{-\infty}^{+\infty}|\o|^{2N}|X(\iomega)|^2 d\o <+\infty.
$$
\par
Further,  Lemma \ref{lemma} claims  predicability functions from
$\X(T)$ on the finite horizon $T$ only.  At the same time, Theorem
\ref{Th1} states that are uniquely defined on $\R$ by their values
at a semi-infinite interval. It does not mean that Lemma \ref{lemma}
is strengthened here: prediction on a horizon larger than $T$ would
be a ill-posed problem that yet has a unique solution.

\end{document}